# Moments of volumes of lower-dimensional random simplices are not monotone

Benjamin Reichenwallner[1]

**Abstract** In a $d$-dimensional convex body $K$, for $n \leq d+1$, random points $X_0, \ldots, X_{n-1}$ are chosen according to the uniform distribution in $K$. Their convex hull is a random $(n-1)$-simplex with probability 1. We denote its $(n-1)$-dimensional volume by $V_{K[n]}$. The $k$-th moment of the $(n-1)$-dimensional volume of a random $(n-1)$-simplex is monotone under set inclusion, if $K \subseteq L$ implies that the $k$-th moment of $V_{K[n]}$ is not larger than that of $V_{L[n]}$. Extending work of Rademacher (Mathematika 58:77–91, 2012) and Reichenwallner and Reitzner (Mathematika 62:949–958, 2016), it is shown that for $n \leq d$, the moments of $V_{K[n]}$ are not monotone under set inclusion. As a consequence, the nonmonotonicity of the expected surface area of the convex hull of $n \leq d+1$ uniform random points in a $d$-dimensional convex body follows.

**Mathematics Subject Classification** Primary 60D05 · Secondary 52A22

## 1 Introduction

For a $d$-dimensional convex body $K$ and $n \leq d+1$, we denote the $(n-1)$-dimensional volume of the convex hull of $n$ independently and uniformly distributed random points $X_0, \ldots, X_{n-1}$ in $K$ by $V_{K[n]}$. Since the $n$ points are not contained in any $(n-2)$- dimensional hyperplane with probability 1, their convex hull is almost surely an $(n-1)$-simplex. In 2006, Meckes [5] asked whether for a pair of convex bodies $K, L \subseteq \mathbb{R}^d$, $K \subseteq L$ would imply

$$\mathbb{E} V_{K[d+1]} \leq \mathbb{E} V_{L[d+1]}. \tag{1}$$

Interest in this question comes from the fact that it would imply a positive solution to the slicing problem, see e.g. [7]. For a more general statement of the conjecture we refer to [9].

In 2012, Rademacher [7] showed that Meckes' conjecture (1) is not true in general, since there exist counterexamples for each dimension $d \geq 4$. His results were generalized to higher moments of $V_{K[d+1]}$ by Reichenwallner and Reitzner [8], and the expected value in dimension three was handled in [4]:

**Theorem 1** *Let $k \in \mathbb{N}$ and $d \geq 2$. For two-dimensional convex bodies $K$ and $L$, $K \subseteq L$ implies that $\mathbb{E} V_{K[d+1]} \leq \mathbb{E} V_{L[d+1]}$ and $\mathbb{E} V^2_{K[d+1]} \leq \mathbb{E} V^2_{L[d+1]}$. Unless $d = 2$ and $k \in \{1, 2\}$, there exist two $d$-dimensional convex bodies $K, L$ satisfying $K \subseteq L$ and $\mathbb{E} V^k_{K[d+1]} > \mathbb{E} V^k_{L[d+1]}$.*

Note that the problem is trivial in dimension one. In this paper, we perform another extension by considering convex hulls of $n \leq d+1$ random points in $d$-dimensional convex bodies. In order to keep the statement of the result simple, we restrict to $n \leq d$ and refer to the theorem above for the case $n = d+1$.

**Theorem 2** *Let $d, n, k \in \mathbb{N}$ with $d \geq 2$ and $2 \leq n \leq d$. Then there exist two $d$-dimensional convex bodies $K, L$ satisfying $K \subseteq L$ and $\mathbb{E} V^k_{K[n]} > \mathbb{E} V^k_{L[n]}$.*

Due to Rademacher [7], it is already known that the moments of $V_{K[n]}$ are not monotone under set inclusion in general. However, it might still come as a surprise that for $n \leq d$, there exist counterexamples for *any* moment in *any* dimension, particularly in dimension two, where the first two moments of the area of a random triangle are monotone.

Theorem 1 is shown using an equivalence statement, first stated in Lemma 12 in [7]. In consideration of Proposition 16 ibidem, it is clear that the lemma can be generalized in the following way:

**Lemma 1** *For $d, n, k \in \mathbb{N}$ with $2 \leq n \leq d+1$, monotonicity under inclusion of the map $K \mapsto \mathbb{E} V^k_{K[n]}$, where $K$ ranges over all $d$-dimensional convex bodies, holds if and only if we have for each convex body $K \subseteq \mathbb{R}^d$ and for each $x \in \mathrm{bd}\, K$ that $\mathbb{E} V^k_{K[n]} \leq \mathbb{E} V^k_{K[n],x}$.*

Here, $V_{K[n],x}$ denotes the $(n-1)$-dimensional volume of a random $(n-1)$-simplex, which is the convex hull of a fixed point $x$ and $n-1$ independently and uniformly distributed random points in $K$. The lemma allows us to consider one convex body $K$, rather than a pair of convex bodies, and compute two different moments: the moment of the $(n-1)$-dimensional volume of a random $(n-1)$-simplex in $K$ as well as the same, but fixing one of the $n$ points to be a point on the boundary of $K$, denoted by $\mathrm{bd}\, K$.

For the case $n = d + 1$, Rademacher chooses $K$ to be a $d$-dimensional halfball and $x$ the midpoint of the base which is a $(d - 1)$-dimensional ball. This pair forms the counterexample in Theorem 1 to the monotonicity of all moments for $d \geq 4$, and to the monotonicity of the $k$-th moment for $k \geq 11$ in dimension two and for $k \geq 2$ in dimension three (see [8]; in fact, in dimension two, it does so for $k \geq 8$). Additional counterexamples are given by the triangle, together with the midpoint of one of its edges, for $d = 2$ and $k \geq 3$ in [8] and the tetrahedron, together with the centroid of one of its facets, for $d = 3$ and $k = 1$ in [4].

For our generalization, we use the following lemma, which allows us to create counterexamples inductively:

**Lemma 2** *Let $d, n, k \in \mathbb{N}$ with $2 \leq n \leq d + 1$. If there exist two $d$-dimensional convex bodies $K, L$ satisfying $K \subseteq L$ and $\mathbb{E}V_{K[n]}^k > \mathbb{E}V_{L[n]}^k$, there also exist two $(d+1)$-dimensional convex bodies $K', L'$ satisfying $K' \subseteq L'$ and $\mathbb{E}V_{K'[n]}^k > \mathbb{E}V_{L'[n]}^k$.*

We use the counterexamples mentioned above as well as Lemma 2 to get results on the monotonicity of the moments of $V_{K[n]}$, where $K$ ranges over all $d$-dimensional convex bodies, $d \geq 2$ and $2 \leq n \leq d$. Since, obviously, monotonicity holds in dimension 1, in order to show that the moments of the length of a random chord in $d$-dimensional convex bodies are not monotone for $d > 1$, we create another counterexample by computing the moments of the length of a random chord inside a rectangular triangle where both points are randomly chosen as well as the same where one point is fixed to be the midpoint of the hypotenuse.

**Proposition 1** *For the triangle $T_2$ with vertices $(0, 0)$, $(1, 0)$ and $(0, 1)$, and $c_2 = (1/2, 1/2)$ the midpoint of the edge $\{(x, y) : x, y \geq 0, x + y = 1\}$, we have that $\mathbb{E}V_{T_2[2], c_2}^k < \mathbb{E}V_{T_2[2]}^k$ for each $k \in \mathbb{N}$.*

Concerning the case $n = 3, k = 1$, we need to construct a pair of three-dimensional convex bodies, or one three-dimensional convex body and a suitable point on its boundary, leading to a counterexample. Already the computation of the expected volume of a random tetrahedron in dimension three is a hard problem, and the same is supposed to be true for the expected area of a random triangle in dimension three.

The only three-dimensional convex sets where the expected volume of a random full-dimensional simplex is known are the ball [6], the cube [12] and the tetrahedron [1], and the ball is the only one for which the expected area of a random triangle inside is known. Since numerical simulations suggest that the ball cannot form a counterexample in our concerns, and neither can the cube, we consider random triangles in a tetrahedron, which nothing is known about. A suitable choice for $x$ should be the center $c$ of one of the facets, but we have to take into account that the considered moments depend on the shape of $T$ as well as of the facet which $c$ lies in. Additionally, already the determination of the expected volume of a random tetrahedron in a tetrahedron $T$, done by Buchta and Reitzner [1], was extremely hard, and it seems that a computation of $\mathbb{E}V_{T[3]}$ and $\mathbb{E}V_{T[3], c}$ is out of reach. Neverless we will prove the following proposition.

**Proposition 2** *For the tetrahedron $T_3$ with vertices $(0, 0, 0)$, $(1, 0, 0)$, $(0, 1, 0)$ and $(0, 0, 1)$, and $c_3 = (1/3, 1/3, 1/3)$ the centroid of the facet $\{(x, y, z) : x, y, z \geq 0, x + y + z = 1\}$, we have that $\mathbb{E}V_{T_3[3], c_3} < \mathbb{E}V_{T_3[3]}$.*

A combination of this result and a direct computation of $\mathbb{E}V^2_{T_3[3]}$ and $\mathbb{E}V^2_{T_3[3],c_3}$ with Rademacher's Lemma 1 yields Theorem 2. Proposition 2 is obtained by a combination of methods from stochastic geometry with results from approximation theory. In the background, we use a result about the precise approximation of the square root function by suitable polynomials.

A next step might naturally be given by a generalization of Theorems 1 and 2 to intrinsic volumes, among them the surface area. While a derivation of results on its higher moments seems to be rather involved, it appears that the nonmonotonicity of the expected surface area of a random simplex, eventually not of full dimension, is a direct consequence of Theorem 2. For a $d$-dimensional convex body $K$ and $2 \leq n \leq d+1$, let $S_{K[n]}$ denote the surface area of the convex hull of $n$ points that are chosen according to uniform distribution in $K$. Since the (relative) boundary of an $(n-1)$-simplex is given by the union of $n$ facets and the expected value is additive, the following corollary holds true:

**Corollary 1** *Let $d, n \in \mathbb{N}$ with $d \geq 2$ and $2 \leq n \leq d+1$. Then there exist two $d$-dimensional convex bodies $K, L$ satisfying $K \subseteq L$ and $\mathbb{E}S_{K[n]} > \mathbb{E}S_{L[n]}$.*

More generally, it shall be mentioned that the expected value of the sum of the $m$-dimensional volumes of all $m$-faces of the convex hull of $n \geq m+2$ uniform random points in a convex body $K$ of dimension $d \geq n-1$ is not monotone under set inclusion.

Finally, we consider a similar problem: Assume that $n$ points are uniformly distributed on the boundary of a convex body $K$ and let $\mathcal{V}_{K[n]}$ be the $(n-1)$-dimensional volume of their convex hull. Then we can use a technique similar to that in the background of Lemma 2 to show that, in general, monotonicity under set inclusion does not hold for the moments of $(n-1)$-dimensional volumes of this type of random $(n-1)$-simplices either. In detail, we show the following:

**Theorem 3** *Let $d, n, k \in \mathbb{N}$ with $d \geq 3$ and $2 \leq n \leq d$. Unless $d = 3$, $n = 3$ and $k \in \{1, 2\}$, there exist two $d$-dimensional convex bodies $K, L$ satisfying $K \subseteq L$ and $\mathbb{E}\mathcal{V}^k_{K[n]} > \mathbb{E}\mathcal{V}^k_{L[n]}$.*

This paper is organized in the following way. We prove Lemma 2 in Sect. 2. In Sect. 3, we compute the moments of the lengths of two types of random chords inside a triangle as described above, and thereafter we give bounds for the expected areas of random triangles in a tetrahedron. These results will be used in Sect. 5 for the proof of Theorem 2. Finally, Theorem 3 will be proved in Sect. 6.

As a general reference for the tools and results we need in the following, we refer to the book on Stochastic and Integral Geometry by Schneider and Weil [10]. More recent surveys on random polytopes are due to Hug [3] and Reitzner [9].

## 2 Proof of Lemma 2

Let $d, n, k \in \mathbb{N}$ with $2 \leq n \leq d+1$ and assume that there exist two $d$-dimensional convex bodies $K$ and $L$ satisfying $K \subseteq L$ and $\mathbb{E}V^k_{K[n]} > \mathbb{E}V^k_{L[n]}$. Our goal is to construct two $(d+1)$-dimensional convex bodies $K'$, $L'$ satisfying similar conditions. For $\varepsilon > 0$, we define two $(d+1)$-dimensional convex bodies by

$$K_\varepsilon := K \times [0, \varepsilon], \quad L_\varepsilon := L \times [0, \varepsilon].$$

Then, $K_\varepsilon \subseteq L_\varepsilon$ for all $\varepsilon$, and if $\mathbb{E}V^k_{K_\varepsilon[n]}$ and $\mathbb{E}V^k_{L_\varepsilon[n]}$ converge to $\mathbb{E}V^k_{K[n]}$ and $\mathbb{E}V^k_{L[n]}$, resp., for $\varepsilon$ tending to zero, there exists an $\varepsilon_0 > 0$ such that $\mathbb{E}V^k_{K_{\varepsilon_0}[n]} > \mathbb{E}V^k_{L_{\varepsilon_0}[n]}$.

It remains to proof the continuity of the map $\varepsilon \mapsto \mathbb{E}V^k_{K_\varepsilon[n]}$ in $\varepsilon = 0$. We express the $(n-1)$-dimensional volume of an $(n-1)$-simplex with vertices $x_0, \ldots, x_{n-1}$ using a matrix that contains the vectors spanning the $(n-1)$-simplex in its columns:

$$\mathrm{vol}_{n-1} \mathrm{conv}(x_0, x_1, \ldots, x_{n-1})$$
$$= \frac{1}{(n-1)!}\sqrt{\det\left((x_1 - x_0, \ldots, x_{n-1} - x_0)^T (x_1 - x_0, \ldots, x_{n-1} - x_0)\right)}.$$

Then we have for $\varepsilon > 0$:

$$\mathbb{E}V^k_{K_\varepsilon[n]} = \frac{1}{\mathrm{vol}\, K_\varepsilon^n} \int_{K_\varepsilon^n} \mathrm{vol}_{n-1} \mathrm{conv}(x_0, \ldots, x_{n-1})^k \, d(x_0, \ldots, x_{n-1})$$
$$= \frac{1}{\varepsilon^n \mathrm{vol}_d K^n} \frac{1}{((n-1)!)^k} \int_{K^n} \int_{[0,\varepsilon]^n} \det B^{k/2} \, d(z_0, \ldots, z_{n-1}) \, d(x_0, \ldots, x_{n-1}),$$

where $B = (b_{ij})_{1 \le i, j \le n-1}$ with

$$b_{ij} = \sum_{m=1}^{d} \left(x_i^{(m)} - x_0^{(m)}\right)\left(x_j^{(m)} - x_0^{(m)}\right) + (z_i - z_0)(z_j - z_0),$$

writing $x_i = (x_i^{(1)}, \ldots, x_i^{(d)}, z_i)$. Now, fix the first $d$ coordinates of the points $x_0, \ldots, x_{n-1}$ and let $\zeta := \min_{z_i} \det B$ and $Z := \max_{z_i} \det B$, where the minimum and the maximum are taken over all possible choices of $(z_0, \ldots, z_{n-1}) \in [0, \varepsilon]^n$. The extrema exist since the determinant is continuous and the interval $[0, \varepsilon]$ is compact. We have:

$$\frac{1}{\mathrm{vol}_d K^n} \frac{1}{((n-1)!)^k} \int_{K^n} \zeta^{k/2} \, d(x_0, \ldots, x_{n-1}) \le \mathbb{E}V^k_{K_\varepsilon[n]}$$
$$\le \frac{1}{\mathrm{vol}_d K^n} \frac{1}{((n-1)!)^k} \int_{K^n} Z^{k/2} \, d(x_0, \ldots, x_{n-1}).$$

If $\varepsilon$ tends to zero, the coordinates $z_i$ do so as well, and all summands of $\zeta$ and $Z$ containing a factor $z_i$ vanish in the limit. Therefore, we get that

$$\mathbb{E}V^k_{K[n]} \le \lim_{\varepsilon \to 0} \mathbb{E}V^k_{K_\varepsilon[n]} \le \mathbb{E}V^k_{K[n]},$$

which implies that $\lim_{\varepsilon \to 0} \mathbb{E}V^k_{K_\varepsilon[n]} = \mathbb{E}V^k_{K[n]}$. □

## 3 Random chords in a triangle

Essential for our investigations of the case $n = 2$ are the moments of two different types of random distances in a given triangle $T$. In detail, we need the moments of the length of a random chord in a right-angled triangle, which arises as the segment between two uniformly chosen random points in this triangle, as well as the moments of the distance of a random point in this triangle from the midpoint of its hypotenuse, which is equivalent to the length of a random chord where one point is fixed to be this midpoint.

### 3.1 Moments of the distance of a random point from a point on the boundary

We first prove a lemma that will be useful for the computation of the moments of the distance of a random point from a fixed point on the boundary, since we can dissect a triangle by a line through this fixed point.

**Lemma 3** *Let $T$ be a triangle with vertices $A$, $B$ and $C$. Denote the lengths of the edges opposite to $A$, $B$ and $C$ by $a$, $b$ and $c$, resp., and the angles in these corners by $\alpha$, $\beta$ and $\gamma$, resp. Then it holds for the moments of the distance of a random point in $T$ from $A$:*

$$\mathbb{E} V^k_{T[2],A} = \frac{2(c \sin \beta)^{k+1}}{(k+2)a} \left( \cos \beta \, {}_2F_1 \left( \frac{1}{2}, \frac{k+3}{2}; \frac{3}{2}; \cos \beta^2 \right) \right.$$
$$\left. - \cos(\alpha + \beta) \, {}_2F_1 \left( \frac{1}{2}, \frac{k+3}{2}; \frac{3}{2}; \cos(\alpha + \beta)^2 \right) \right).$$

Here, the function

$$ {}_2F_1(a_1, a_2; b; x) = \sum_{k=0}^{\infty} \frac{(a_1)_k (a_2)_k}{(b)_k} \frac{x^k}{k!} $$

with $(a)_k = a(a+1) \cdots (a+k-1)$ is a hypergeometric function.

*Proof* Since the expected length of a random chord in a triangle does not depend on the position of the latter, we build our coordinate system ensuring that $A$ is identical to the origin $o$ and the edge $c$ lies on the $x$-axis. Then:

$$\mathbb{E} V^k_{T[2],o} = \frac{1}{\text{vol } T} \int_T ||x||^k \, dx.$$

Note that vol $T = ac \sin \beta /2$. Using polar coordinates, we get that

$$\mathbb{E} V^k_{T[2],o} = \frac{2}{ac \sin \beta} \int_0^\alpha \int_0^{l(\varphi)} r^{k+1} \, dr \, d\varphi = \frac{1}{k+2} \frac{2}{ac \sin \beta} \int_0^\alpha l(\varphi)^{k+2} \, d\varphi,$$

where $l(\varphi)$ denotes the length of the intersection of the triangle with the line whose angle to the $x$-axis is $\varphi$. By the sine rules, $l(\varphi) = c \sin \beta / \sin(\beta + \varphi)$, and

$$\mathbb{E} V_{T[2],o}^k = \frac{2(c \sin \beta)^{k+1}}{(k+2)a} \int_0^\alpha \frac{1}{\sin(\beta + \varphi)^{k+2}} \, d\varphi.$$

With

$$\int \frac{1}{\sin \varphi^{k+2}} \, d\varphi = -\cos \varphi \, _2F_1\left(\frac{1}{2}, \frac{k+3}{2}; \frac{3}{2}; \cos \varphi^2\right),$$

we get that

$$\mathbb{E} V_{T[2],o}^k = \frac{2(c \sin \beta)^{k+1}}{(k+2)a} \left( \cos \beta \, _2F_1\left(\frac{1}{2}, \frac{k+3}{2}; \frac{3}{2}; \cos \beta^2\right) \right.$$
$$\left. - \cos(\alpha + \beta) \, _2F_1\left(\frac{1}{2}, \frac{k+3}{2}; \frac{3}{2}; \cos(\alpha + \beta)^2\right) \right).$$

□

Given this lemma, we can compute the moments of the distance of a random point from a point on the boundary of a triangle.

**Proposition 3** *Let $T$ be a triangle with vertices $A$, $B$ and $C$. Denote the lengths of the edges opposite to $A$, $B$ and $C$ by $a$, $b$ and $c$, resp., and the angles in these corners by $\alpha$, $\beta$ and $\gamma$, resp. Furthermore, let $D$ be a point lying on the edge between $A$ and $B$ and $c_1$ the distance of $D$ from $A$. Then it holds for the moments of the distance of a random point in $T$ from $D$:*

$$\mathbb{E} V_{T[2],D}^k = \frac{2}{(k+2)bc \sin \alpha} \left( (c_1 \sin \alpha)^{k+2} \cos \alpha \, _2F_1\left(\frac{1}{2}, \frac{k+3}{2}; \frac{3}{2}; \cos \alpha^2\right) \right.$$
$$- (c_1 \sin \alpha)^{k+2} \cos(\alpha + \delta) \, _2F_1\left(\frac{1}{2}, \frac{k+3}{2}; \frac{3}{2}; \cos(\alpha + \delta)^2\right)$$
$$+ ((c - c_1) \sin \beta)^{k+2} \cos \beta \, _2F_1\left(\frac{1}{2}, \frac{k+3}{2}; \frac{3}{2}; \cos \beta^2\right)$$
$$\left. + ((c - c_1) \sin \beta)^{k+2} \cos(\beta - \delta) \, _2F_1\left(\frac{1}{2}, \frac{k+3}{2}; \frac{3}{2}; \cos(\beta - \delta)^2\right) \right),$$

*where*

$$\delta = \arcsin\left(b \sin \alpha / \sqrt{c_1^2 + b^2 - 2c_1 b \cos \alpha}\right).$$

*Proof* We split $T$ into two parts by the line segment between $C$ and $D$. Denote the length of this line segment by $d$ and the angle between $AD$ and $DC$ by $\delta$. Then, by cosine and sine rule,

$$d = \sqrt{c_1^2 + b^2 - 2c_1 b \cos\alpha} \quad \text{and} \quad \delta = \arcsin(b \sin\alpha/d).$$

If we denote the triangle $ACD$ by $T_1$ and the triangle $BCD$ by $T_2$, we have:

$$\mathbb{E}V_{T[2],D}^k = \frac{\text{vol } T_1}{\text{vol } T} \mathbb{E}V_{T_1[2],D}^k + \frac{\text{vol } T_2}{\text{vol } T} \mathbb{E}V_{T_2[2],D}^k.$$

The previous lemma yields the result. □

In order to find a counterexample to the monotonicity of the moments of the length of random chords, it suffices to consider one specific triangle. We choose a right-angled, isosceles one and apply Proposition 3 to this triangle.

**Corollary 2** *For the triangle $T_2$ with vertices $(0, 0)$, $(0, 1)$ and $(1, 0)$, and $c_2 = (1/2, 1/2)$, we have:*

$$\mathbb{E}V_{T_2[2],c_2}^k = \frac{2^{1/2-k}}{k+2} \, {}_2F_1\left(\frac{1}{2}, \frac{k+3}{2}; \frac{3}{2}; \frac{1}{2}\right).$$

We give some values of $\mathbb{E}V_{T_2[2],c_2}^k$ for this specific triangle:

$$\mathbb{E}V_{T_2[2],c_2} = \frac{1}{6\sqrt{2}}(2 + \sqrt{2} \operatorname{arcsinh} 1) = 0.3825\ldots,$$

$$\mathbb{E}V_{T_2[2],c_2}^2 = \frac{1}{6} = 0.1666\ldots,$$

$$\mathbb{E}V_{T_2[2],c_2}^3 = \frac{1}{160\sqrt{2}}(14 + 3\sqrt{2} \operatorname{arcsinh} 1) = 0.0783\ldots,$$

$$\mathbb{E}V_{T_2[2],c_2}^4 = \frac{7}{180} = 0.0388\ldots$$

### 3.2 Moments of the length of a random chord

In this subsection, we complete our investigations of the triangle by a computation of the moments of the length of a random chord in a triangle, where both endpoints of the chord are uniformly distributed points in the triangle.

**Proposition 4** *Let $T$ be a triangle with vertices $E_1$, $E_2$ and $E_3$. Denote the length of the edge opposite to $E_i$ by $e_i$ and the angle in this corner by $\eta_i$, $i = 1, 2, 3$. Then it holds for the moments of the distance of two random point in $T$:*

$$\mathbb{E}V_{T[2]}^k = \frac{8}{(k+2)(k+3)(k+4)} \frac{1}{(e_1 e_3 \sin\eta_2)^2} \sum_{i \neq j} \sin\eta_i^{k+3} e_j^{k+4}$$

$$\times \left[\cos\eta_i \cos\eta_j \, {}_2F_1\left(\frac{1}{2}, \frac{k+3}{2}; \frac{3}{2}; \cos\eta_j^2\right)\right.$$

$$+ \cos \eta_i^2 \, _2F_1\left(\frac{1}{2}, \frac{k+3}{2}; \frac{3}{2}; \cos \eta_i^2\right)$$

$$+ \frac{\sin \eta_i}{k+2} \left(\frac{1}{\sin \eta_j^{k+2}} - \frac{1}{\sin \eta_i^{k+2}}\right)\Bigg].$$

*Proof* As before, the area of the triangle is given by $e_1 e_3 \sin \eta_2 / 2$. Using the affine Blaschke-Petkantschin formula—see e.g. [10]—we integrate over all lines $H$ intersecting the triangle:

$$\mathbb{E} V_{T[2]}^k = \frac{1}{\text{vol } T^2} \int_{T^2} \|x_0 - x_1\|^k \, d(x_0, x_1)$$

$$= \frac{4}{(e_1 e_3 \sin \eta_2)^2} \int_{A(2,1)} \int_{(H \cap T)^2} \|x_0 - x_1\|^k \, d(x_0, x_1) \, dH,$$

where $A(2, 1)$ denotes the affine Grassmannian of lines in $\mathbb{R}^2$. We represent a line $H$ by its unit normal vector $u \in S_1$ and its distance $t > 0$ from the origin and we therefore denote the line by

$$H_{t,u} = \{x \in \mathbb{R}^2 : \langle x, u \rangle = t\}.$$

We choose the normalization of the Haar measure $dH$ in such a way that $dH = dt \, du$, where $dt$ and $du$ correspond to Lebesgue measures in $\mathbb{R}$ and $S_1$. We write the appearing integral as an expectation to get

$$\mathbb{E} V_{T[2]}^k = \frac{4}{(e_1 e_3 \sin \eta_2)^2} \int_{S_1} \int_0^\infty \int_{(H_{t,u} \cap T)^2} \|x_0 - x_1\|^{k+1} \, d(x_0, x_1) \, dt \, du$$

$$= \frac{4}{(e_1 e_3 \sin \eta_2)^2} \int_{S_1} \int_0^\infty \text{vol}(H_{t,u} \cap T)^2 \, \mathbb{E} V_{H_{t,u} \cap T}^{k+1} \, dt \, du.$$

The $(k + 1)$-st moment of the distance of two random points in an interval $I$ of length $l$ has already been computed by Solomon [11], amongst others, to be given by

$$\mathbb{E} V_I^{k+1} = \frac{2 \, l^{k+1}}{(k+2)(k+3)}.$$

Hence, we obtain

$$\mathbb{E} V_{T[2]}^k = \frac{8}{(k+2)(k+3)} \frac{1}{(e_1 e_3 \sin \eta_2)^2} \int_{S_1} \int_0^\infty \text{vol}(H_{t,u} \cap T)^{k+3} \, dt \, du.$$

A line $H_{t,u}$ that intersects the triangle $T$ a.s. meets exactly two edges of $T$. It splits $T$ into a triangle and a quadrangle. We say that $H_{t,u}$ cuts off a vertex $E_i$ from $T$ if $E_i$ is contained in the triangular part. Furthermore, we write

$$\mathcal{I}^{(i)} = \int_{S_1} \int_0^\infty \mathbb{1}(H_{t,u} \text{ cuts off } E_i \text{ from } T) \operatorname{vol}(H_{t,u} \cap T)^{k+3} \, dt \, du$$

for $i = 1, 2, 3$, which gives

$$\mathbb{E} V_{T[2]}^k = \frac{8}{(k+2)(k+3)} \frac{1}{(e_1 e_3 \sin \eta_2)^2} \left( \mathcal{I}^{(1)} + \mathcal{I}^{(2)} + \mathcal{I}^{(3)} \right).$$

We state the following lemma, which will be proved right after the end of the proof of this proposition:

**Lemma 4** *It holds for $i = 1, 2, 3$:*

$$\mathcal{I}^{(i)} = \frac{\sin \eta_i^{k+3}}{k+4} \left( e_j^{k+4} \left[ \cos \eta_i \cos \eta_j \, {}_2F_1\left(\frac{1}{2}, \frac{k+3}{2}; \frac{3}{2}; \cos \eta_j^2\right) \right. \right.$$

$$+ \cos \eta_i^2 \, {}_2F_1\left(\frac{1}{2}, \frac{k+3}{2}; \frac{3}{2}; \cos \eta_i^2\right)$$

$$\left. + \frac{\sin \eta_i}{k+2} \left( \frac{1}{\sin \eta_j^{k+2}} - \frac{1}{\sin \eta_i^{k+2}} \right) \right]$$

$$+ e_l^{k+4} \left[ \cos \eta_i \cos \eta_l \, {}_2F_1\left(\frac{1}{2}, \frac{k+3}{2}; \frac{3}{2}; \cos \eta_l^2\right) \right.$$

$$+ \cos \eta_i^2 \, {}_2F_1\left(\frac{1}{2}, \frac{k+3}{2}; \frac{3}{2}; \cos \eta_i^2\right)$$

$$\left. \left. + \frac{\sin \eta_i}{k+2} \left( \frac{1}{\sin \eta_l^{k+2}} - \frac{1}{\sin \eta_i^{k+2}} \right) \right] \right).$$

This lemma immediately yields the result. □

*Proof of Lemma 4* Let $i \in \{1, 2, 3\}$. We substitute $u = (\cos \varphi, \sin \varphi)$ such that the line $\varphi = 0$ contains the edge between $E_i$ and $E_j$ for $j \neq i$ fixed, and by the well-known fact that $du = d\varphi$ and with $H_{t,\varphi} = \{(x, y) \in \mathbb{R}^2 : x \cos \varphi + y \sin \varphi = t\}$, we get the integral

$$\mathcal{I}^{(i)} = \int_{-\pi}^{\pi} \int_0^\infty \mathbb{1}(H_{t,\varphi} \text{ cuts off } E_i \text{ from } T) \operatorname{vol}(H_{t,\varphi} \cap T)^{k+3} \, dt \, d\varphi.$$

Let $l \in \{1, 2, 3\} \setminus \{i, j\}$. A line $H_{t,\varphi}$ cuts off $E_i$ from $T$ if and only if $\varphi \in [\eta_i - \pi/2, \pi/2]$ and

(i) $0 \leq t \leq e_j \cos(\eta_1 - \varphi)$, if $\varphi \leq \frac{\pi}{2} - \eta_j$,
(ii) $0 \leq t \leq e_l \cos \varphi$, if $\varphi \geq \frac{\pi}{2} - \eta_j$.

We use the notation $l(t, \varphi) = \mathrm{vol}(H_{t,\varphi} \cap T)$ for the length of the intersection of $T$ with a line $H_{t,\varphi}$. Basic trigonometric considerations show that

$$l(t, \varphi) = t(\tan(\eta_i - \varphi) + \tan \varphi) = t \frac{\sin \eta_i}{\cos \varphi \cos(\eta_i - \varphi)}.$$

This leads to the expression

$$\mathcal{I}^{(i)} = \int_{\eta_i - \pi/2}^{\pi/2 - \eta_j} \left( \frac{\sin \eta_i}{\cos \varphi \cos(\eta_i - \varphi)} \right)^{k+3} \int_0^{e_j \cos(\eta_i - \varphi)} t^{k+3} \, dt \, d\varphi$$

$$+ \int_{\pi/2 - \eta_j}^{\pi/2} \left( \frac{\sin \eta_i}{\cos \varphi \cos(\eta_i - \varphi)} \right)^{k+3} \int_0^{e_l \cos \varphi} t^{k+3} \, dt \, d\varphi$$

$$= \frac{\sin \eta_i^{k+3}}{k+4} \left( e_j^{k+4} \int_{\eta_i - \pi/2}^{\pi/2 - \eta_j} \frac{\cos(\eta_i - \varphi)}{\cos \varphi^{k+3}} \, d\varphi + e_l^{k+4} \int_{\pi/2 - \eta_j}^{\pi/2} \frac{\cos \varphi}{\cos(\eta_i - \varphi)^{k+3}} \, d\varphi \right).$$

With a substitution $\psi = \eta_i - \varphi$ and with $\cos(\eta_i - \varphi) = \cos \eta_i \cos \varphi + \sin \eta_i \sin \varphi$, we can compute the outer integral:

$$\mathcal{I}^{(i)} = \frac{\sin \eta_i^{k+3}}{k+4} \left( e_j^{k+4} \left[ \cos \eta_i \cos \eta_j \, {}_2F_1 \left( \frac{1}{2}, \frac{k+3}{2}; \frac{3}{2}; \cos \eta_j^2 \right) \right. \right.$$
$$+ \cos \eta_i^2 \, {}_2F_1 \left( \frac{1}{2}, \frac{k+3}{2}; \frac{3}{2}; \cos \eta_i^2 \right)$$
$$\left. + \frac{\sin \eta_i}{k+2} \left( \frac{1}{\sin \eta_j^{k+2}} - \frac{1}{\sin \eta_i^{k+2}} \right) \right]$$
$$+ e_l^{k+4} \left[ \cos \eta_i \cos \eta_l \, {}_2F_1 \left( \frac{1}{2}, \frac{k+3}{2}; \frac{3}{2}; \cos \eta_l^2 \right) \right.$$
$$+ \cos \eta_i^2 \, {}_2F_1 \left( \frac{1}{2}, \frac{k+3}{2}; \frac{3}{2}; \cos \eta_i^2 \right)$$
$$\left. \left. + \frac{\sin \eta_i}{k+2} \left( \frac{1}{\sin \eta_l^{k+2}} - \frac{1}{\sin \eta_i^{k+2}} \right) \right] \right).$$

$\square$

Again, we apply this proposition to the right-angled, isosceles triangle.

**Corollary 3** *For the triangle $T_2$ with vertices $(0, 0)$, $(0, 1)$ and $(1, 0)$, we have:*

$$\mathbb{E}V_{T_2[2]}^k = \frac{2^{(5-k)/2} + 2^{7/2}}{(k+2)(k+3)(k+4)} \, {}_2F_1\left(\frac{1}{2}, \frac{k+3}{2}; \frac{3}{2}; \frac{1}{2}\right).$$

As before, we give some values of $\mathbb{E}V_{T_2[2]}^k$ for this triangle:

$$\mathbb{E}V_{T_2[2]} = \frac{1+2\sqrt{2}}{30}(2 + \sqrt{2}\operatorname{arcsinh} 1) = 0.4142\ldots,$$

$$\mathbb{E}V_{T_2[2]}^2 = \frac{2}{9} = 0.2222\ldots,$$

$$\mathbb{E}V_{T_2[2]}^3 = \frac{1+4\sqrt{2}}{840}(14 + 3\sqrt{2}\operatorname{arcsinh} 1) = 0.1405\ldots,$$

$$\mathbb{E}V_{T_2[2]}^4 = \frac{1}{10} = 0.1.$$

## 4 Random triangles in a tetrahedron

In this section, our aim is to show that for a specific tetrahedron $T$ and the center $c$ of one of its facets, we have

$$\mathbb{E}V_{T[3]} > \mathbb{E}V_{T[3],c}. \tag{2}$$

Therefore, we want to find polynomials $P(x) = \sum_{i=0}^n a_i x^i$ and $Q(x) = \sum_{i=0}^l b_i x^i$ such that $P(x) \leq |x|$ and $Q(x) \geq |x|$ for relevant ranges of values $x \in \mathbb{R}$ and consequently

$$\mathbb{E}V_{T[3]} \geq \mathbb{E}P(V_{T[3]}) > \mathbb{E}Q(V_{T[3],c}) \geq \mathbb{E}V_{T[3],c},$$

which suffices in order to show that the pair $(T, c)$ is a counterexample to monotonicity of the expected area of a random triangle in dimension three. Simulations suggest that

$$\mathbb{E}V_{T[3]} = 0.0592\ldots \text{ and } \mathbb{E}V_{T[3],c} = 0.0466\ldots \tag{3}$$

and conjecture (2) should hold true.

### 4.1 A lower bound for the expected area of a random triangle in a tetrahedron

Let $T$ be a tetrahedron. For random points $X_0, X_1, X_2 \in T$, we write $X_i = (x_i, y_i, z_i)$. The volume of the triangle with vertices $X_0$, $X_1$ and $X_2$ is given by

$$V_{T[3]} = \frac{1}{2}\sqrt{\det\left(\begin{pmatrix} x_1 - x_0 & y_1 - y_0 & z_1 - z_0 \\ x_2 - x_0 & y_2 - y_0 & z_2 - z_0 \end{pmatrix} \begin{pmatrix} x_1 - x_0 & x_2 - x_0 \\ y_1 - y_0 & y_2 - y_0 \\ z_1 - z_0 & z_2 - z_0 \end{pmatrix}\right)}$$

$$= \frac{1}{2} D(x_0, \ldots, z_2)^{1/2},$$

and this is the square root of a polynomial of degree four in the coordinates of $X_0$, $X_1$ and $X_2$. For even moments of $V_{T[3]}$, we get rid of the square root:

$$\mathbb{E} V_{T[3]}^{2k} = \text{vol } T^{-3} 2^{-2k} \int_T \int_T \int_T D(x_0, \ldots, z_2)^k \, d(x_0, y_0, z_0) \, d(x_1, y_1, z_1) \, d(x_2, y_2, z_2). \tag{4}$$

Note that the expectation $\mathbb{E} V_{T[3]}$ depends on the shape of the tetrahedron $T$. We consider the specific tetrahedron

$$T_3 = \{(x, y, z) \in \mathbb{R}^3 : x, y, z \geq 0, x + y + z \leq 1\},$$

i.e., that with vertices $(0, 0, 0)$, $(1, 0, 0)$, $(0, 1, 0)$ and $(0, 0, 1)$. Expanding the determinant, the polynomial $D$ can be written as

$$D(x_0, \ldots, z_2) = \sum_{i=1}^{63} s_i \text{ with } s_i = \varepsilon_i \, 2^{c_i} \, x_0^{e_{i0}} \cdots z_2^{e_{i8}},$$

where $\varepsilon_i \in \{-1, 1\}$, $c_i \in \{0, 1\}$ and $e_{ij} \geq 0$, $j = 0, \ldots, 8$, are explicitly given constants and $\sum_j e_{ij} = 4$ for $i = 1, \ldots, 63$. By the Multinomial Theorem, and using the multinomial coefficient

$$\binom{k}{k_1, \ldots, k_{63}} = \frac{k!}{k_1! \cdots k_{63}!},$$

the $k$-th power of it can be rewritten as

$$D(x_0, \ldots, z_2)^k = \sum_{\sum_1^{63} k_i = k} \binom{k}{k_1, \ldots, k_{63}} \prod_{i=1}^{63} s_i^{k_i}$$

$$= \sum_{\sum_1^{63} k_i = k} (-1)^{k'} 2^{k''} \binom{k}{k_1, \ldots, k_{63}} \prod_{i=0}^{2} x_i^{l_1} y_i^{m_1} z_i^{n_1}. \tag{5}$$

For completeness, we list the exponents appearing in (5):

$k' = k_{19} + \cdots + k_{37} + k_{40} + k_{44} + k_{47} + k_{51} + k_{54} + k_{55} + k_{59} + k_{63},$
$k'' = k_{19} + \cdots + k_{63},$
$l_0 = 2k_1 + \cdots + 2k_4 + 2k_{19} + 2k_{20} + k_{27} + k_{29} + k_{33} + k_{35} + k_{37} + \cdots + k_{48},$
$m_0 = 2k_5 + 2k_9 + 2k_{13} + 2k_{14} + k_{21} + k_{23} + 2k_{25} + 2k_{26} + k_{34} + k_{36} + \cdots + k_{38}$
$\quad + k_{43} + k_{44} + k_{49} + k_{50} + k_{55} + \cdots + k_{60},$
$n_0 = 2k_7 + 2k_{11} + 2k_{15} + 2k_{17} + k_{22} + k_{24} + k_{28} + k_{30} + 2k_{31} + 2k_{32} + k_{40} + k_{41}$
$\quad + k_{46} + k_{47} + k_{52} + k_{53} + k_{55} + k_{56} + k_{58} + k_{59} + k_{61} + k_{62},$

$$l_1 = 2k_5 + \cdots + 2k_8 + 2k_{21} + 2k_{22} + k_{25} + k_{29} + k_{31} + k_{35} + k_{37} + \cdots + k_{42}$$
$$+ k_{49} + \cdots + k_{54},$$
$$m_1 = 2k_1 + 2k_{10} + 2k_{15} + 2k_{16} + k_{19} + k_{23} + 2k_{27} + 2k_{28} + k_{32} + k_{36} + k_{37} + k_{39}$$
$$+ k_{43} + k_{45} + k_{49} + k_{51} + k_{55} + \cdots + k_{57} + k_{61} + \cdots + k_{63},$$
$$n_1 = 2k_3 + 2k_{12} + 2k_{13} + 2k_{18} + k_{20} + k_{24} + k_{26} + k_{30} + 2k_{33} + 2k_{34} + k_{40} + k_{42}$$
$$+ k_{46} + k_{48} + k_{52} + k_{54} + k_{55} + k_{57} + k_{58} + k_{60} + \cdots + k_{63},$$
$$l_2 = 2k_9 + 2k_{10} + 2k_{11} + 2k_{12} + 2k_{23} + 2k_{24} + k_{25} + k_{27} + k_{31} + k_{33} + k_{43}$$
$$+ \cdots + k_{54},$$
$$m_2 = 2k_2 + 2k_6 + 2k_{17} + 2k_{18} + k_{19} + k_{21} + 2k_{29} + 2k_{30} + k_{32} + k_{34} + k_{38} + k_{39}$$
$$+ k_{44} + k_{45} + k_{50} + k_{51} + k_{58} + \cdots + k_{63},$$
$$n_2 = 2k_4 + 2k_8 + 2k_{14} + 2k_{16} + k_{20} + k_{22} + k_{26} + k_{28} + 2k_{35} + 2k_{36} + k_{41} + k_{42}$$
$$+ k_{47} + k_{48} + k_{53} + k_{54} + k_{56} + k_{57} + k_{59} + k_{60} + k_{62} + k_{63}.$$

Integration of the monomials over the tetrahedron $T_3$ gives with the substitution $z = t$, $y = s(1-t)$, $x = r(1-s)(1-t)$ that

$$\int_{T_3} x^{l_i} y^{m_i} z^{n_i} \, d(x,y,z) = \int_0^1 r^{l_i} \, dr \int_0^1 s^{m_i} (1-s)^{l_i+1} \, ds \int_0^1 t^{n_i} (1-t)^{l_i+m_i+2} \, dt$$
$$= \frac{1}{l_i+1} B(m_i+1, l_i+2) B(n_i+1, l_i+m_i+3)$$
$$= \frac{l_i! \, m_i! \, n_i!}{(l_i+m_i+n_i+3)!}. \qquad (6)$$

Combining this with Eqs. (4) and (5) gives

$$\mathbb{E} V^{2k}_{T_3[3]} = 6^3 \, 2^{-2k} \sum_{\sum_1^{63} k_i = k} (-1)^{k'} 2^{k''} \binom{k}{k_1, \ldots, k_{63}} \prod_{i=0}^{2} \frac{l_i! \, m_i! \, n_i!}{(l_i+m_i+n_i+3)!}.$$

We list the first five even moments of the area of a random triangle in our tetrahedron $T_3$:

$$\mathbb{E} V^2_{T_3[3]} = \frac{9}{1\,600} = 5.625 \cdot 10^{-3},$$
$$\mathbb{E} V^4_{T_3[3]} = \frac{27}{196\,000} \approx 1.37755 \cdot 10^{-4},$$
$$\mathbb{E} V^6_{T_3[3]} = \frac{3\,161}{379\,330\,560} \approx 8.3331 \cdot 10^{-6},$$
$$\mathbb{E} V^8_{T_3[3]} = \frac{93\,957}{106\,247\,680\,000} \approx 8.8432 \cdot 10^{-7},$$

$$\mathbb{E} V_{T_3[3]}^{10} = \frac{209\,022\,679}{1\,551\,386\,124\,288\,000} \approx 1.34733 \cdot 10^{-7}.$$

We have seen that the $(2k)$-th moment of $V_{T_3[3]}$ can be computed, with fast increasing complexity in $k$. Also note that the area of a triangle in $T_3$ is not larger than $\sqrt{3}/2$, which is the area of the facet $\{(x, y, z) \in \mathbb{R}^3 : x, y, z \geq 0, x + y + z = 1\}$. Hence, we want to approximate the square root function in the interval $[0, 3/4]$ by a polynomial

$$P(x) = \sum_{i=0}^{n} a_i x^i$$

of degree $n$ for some $n \in \mathbb{N}$ such that $P(x) \leq \sqrt{x}$ for all $x \in [0, \frac{3}{4}]$. This way, we can approximate $\mathbb{E} V_{T_3[3]}$ from below by the polynomial $P$:

$$\mathbb{E} V_{T_3[3]} \geq \mathbb{E} P(V_{T_3[3]}).$$

Moreover, the best polynomial for fixed $n \in \mathbb{N}$ can be found via the linear optimization problem

$$\max_{P} \mathbb{E} P(V_{T_3[3]}) = \max_{a_i} \sum_{i=0}^{n} a_i \, \mathbb{E} P(V_{T_3[3]})^{2i} \quad \text{s.t.} \quad P\left(x^2\right) \leq x, \ x \in \left[0, \frac{\sqrt{3}}{2}\right].$$

This constraint is infinite dimensional, but we get an upper bound on $\mathbb{E} P(V_{T_3[3]})$ via the finite dimensional linear program

$$\max_{P} \mathbb{E} P(V_{T_3[3]}) = \max_{a_i} \sum_{i=0}^{n} a_i \, \mathbb{E} V_{T_3[3]}^{2i} \quad \text{s.t.} \quad P\left(x_\ell^2\right) \leq x_\ell, \ x_\ell \in \left[0, \frac{\sqrt{3}}{2}\right],$$
$$\ell = 0, \ldots, L.$$

For $n = 6$ and $L = 200$ equidistant points $x_\ell \in [0, \frac{\sqrt{3}}{2}]$, we numerically compute via Matlab and the optimization toolbox CVX [2] that

$$\mathbb{E} P(V_{T_3[3]}) < 0.04647,$$

and in consideration of the simulations in (3), this yields that we do not get a sufficiently precise estimate using only $n = 6$ even moments.

We solve the above linear program for $n = 7$ and $L = 1\,000$, get the squares of three interpolation nodes with the square root function in $(0, \frac{\sqrt{3}}{2})$ and rationalize these points to

$$\{x_1, x_2, x_3\} = \left\{\frac{2}{19}, \frac{4}{15}, \frac{8}{17}\right\}.$$

Additionally, we take the points $x_0 = 0$ and $x_4 = \frac{47}{54} > \frac{\sqrt{3}}{2}$ as interpolation nodes of our desired polynomial. We want to use these five points to get a polynomial of degree 7. Consider a variation of Lemma 2 in [4] with two conditions

$$P(x_j^2) = x_j \text{ and } P'(x_j^2) = \frac{1}{2x_j}$$

for $j = 1, 2, 3$ and the single conditions $P(x_0^2) = P(0) = 0$ and $P(x_4^2) = x_4$. Then we have with $f(x) = \sqrt{x}$ that

$$f(x) - P(x) = \frac{f^{(n+1)}(\xi)}{(n+1)!} x (x - x_4^2) \prod_{j=1}^{3} (x - x_j^2)^2 \geq 0$$

for some $\xi \in [0, x_4^2]$, since $f^{(n+1)}(x) < 0$, $x \geq 0$ and $x \leq x_4$ for all $x \in [0, \frac{3}{4}]$.

We solve the interpolation problem and finally get a polynomial $P_{\text{cert}}(x) = \sum_{i=1}^{7} a_i x^i$ of degree 7 with explicitly given rational coefficients $a_1, \ldots, a_7$ and the property $P_{\text{cert}} \leq \sqrt{x}$, and this polynomial yields with usage of the even moments computed above that

$$\mathbb{E} V_{T_3[3]} \geq \mathbb{E} P_{\text{cert}}(V_{T_3[3]})$$
$$= \frac{91\,096\,443\,868\,688\,279\,627\,145\,529\,306\,423\,486\,222\,767\,653\,883\,418\,444\,589}{1\,940\,591\,908\,276\,502\,542\,424\,150\,723\,573\,750\,536\,040\,094\,315\,231\,313\,920\,000}$$
$$> 0.046942. \tag{7}$$

### 4.2 An upper bound for the expected area of a random triangle in a tetrahedron where one point is fixed

Now we consider the convex hull of two random points $X_1, X_2$ and one fixed point in a tetrahedron $T$. We write again $X_i = (x_i, y_i, z_i)$ and choose the centroid $c = (x_c, y_c, z_c)$ of one of the facets of $T$ as third point. The volume of the triangle with vertices $X_1, X_2$ and $c$ is given by

$$V_{T[3],c} = \frac{1}{2} \sqrt{\det\left(\begin{pmatrix} x_1 - x_c & y_1 - y_c & z_1 - z_c \\ x_2 - x_c & y_2 - y_c & z_2 - z_c \end{pmatrix} \begin{pmatrix} x_1 - x_c & x_2 - x_c \\ y_1 - y_c & y_2 - y_c \\ z_1 - z_c & z_2 - z_c \end{pmatrix}\right)}$$
$$= \frac{1}{2} D(x_0, \ldots, z_c)^{1/2},$$

and this is again the square root of a polynomial of degree four in the coordinates of $X_1, X_2$ and $c$. As before, we have:

$$\mathbb{E} V_{T[3],c}^{2k} = \text{vol } T^{-2} 2^{-2k} \int_T \int_T D(x_1, \ldots, z_c)^k \, d(x_1, y_1, z_1) \, d(x_2, y_2, z_2). \tag{8}$$

The expectation $\mathbb{E} V_{T[3],c}$ particularly depends on the facet whose centroid we choose to be the fixed vertex. We take again the tetrahedron $T_3$ as defined in Sect. 4.1 and choose $c_3 = (1/3, 1/3, 1/3)$, the centroid of the facet $\{(x, y, z) \in \mathbb{R}^3 : x, y, z \geq 0, x + y + z = 1\}$, to be the fixed vertex of the random triangle.

Expanding the determinant and inserting the coordinates of $c_3$, the polynomial $D$ can be displayed as

$$D(x_1, \ldots, z_{c_3}) = \frac{1}{9} \sum_{i=1}^{54} s_i \quad \text{with} \quad s_i = \varepsilon_i \, 2^{c_i} \, 3^{d_i} \, x_1^{e_{i0}} \cdots z_2^{e_{i5}},$$

where $\varepsilon_i \in \{-1, 1\}$, $c_i, d_i \in \{0, 1, 2\}$ and $e_{ij} \geq 0$, $j = 0, \ldots, 5$, are explicitly given constants and $\sum_j e_{ij} \leq 4$ for $i = 1, \ldots, 54$. Again with usage of the Multinomial Theorem, the series representation of the $k$-th power of it is given by

$$\begin{aligned} D(x_1, \ldots, z_{c_3})^k &= 3^{-2k} \sum_{\sum_1^{54} k_i = k} \binom{k}{k_1, \ldots, k_{54}} \prod_{i=1}^{54} s_i^{k_i} \\ &= 3^{-2k} \sum_{\sum_1^{54} k_i = k} (-1)^{k'} 2^{k''} 3^{k'''} \binom{k}{k_1, \ldots, k_{54}} \prod_{i=1}^{2} x_i^{l_1} y_i^{m_1} z_i^{n_1}, \quad (9) \end{aligned}$$

where we use abbreviations analogously as in the previous proof:

$$\begin{aligned} k' &= k_7 + \cdots + k_{10} + k_{12} + k_{15} + k_{17} + k_{18} + k_{21} + \cdots + k_{33} + k_{52} + \cdots + k_{54}, \\ k'' &= k_1 + \cdots k_6 + 2k_7 + \cdots + 2k_9 + k_{10} + \cdots + k_{45} + k_{52} + \cdots + k_{54}, \\ k''' &= k_{22} + \cdots + k_{45} + 2k_{46} + \cdots + 2k_{54}, \\ l_1 &= 2k_1 + k_7 + k_{10} + \cdots + k_{13} + 2k_{22} + 2k_{23} + k_{28} + k_{32} + k_{34} + \cdots + k_{39} \\ &\quad + 2k_{46} + 2k_{47} + k_{52} + k_{53}, \\ m_1 &= 2k_3 + k_8 + k_{10} + k_{14} + k_{18} + k_{19} + k_{24} + 2k_{26} + 2k_{27} + k_{33} + k_{34} \\ &\quad + k_{38} + k_{40} + k_{42} + \cdots + k_{44} + 2k_{48} + 2k_{50} + k_{52} + k_{54}, \\ n_1 &= 2k_5 + k_9 + k_{12} + k_{16} + k_{18} + k_{20} + k_{25} + k_{29} + 2k_{30} + 2k_{31} + k_{36} \\ &\quad + k_{39} + k_{41} + k_{42} + k_{44} + k_{45} + 2k_{49} + 2k_{51} + k_{53} + k_{54}, \\ l_2 &= 2k_2 + k_7 + k_{14} + \cdots + k_{17} + 2k_{24} + 2k_{25} + k_{26} + k_{30} + k_{34} \\ &\quad + \cdots + k_{37} + k_{40} + k_{41} + 2k_{48} + 2k_{49} + k_{52} + k_{53}, \\ m_2 &= 2k_4 + k_8 + k_{11} + k_{15} + k_{20} + \cdots + k_{22} + 2k_{28} + 2k_{29} + k_{31} + k_{35} \\ &\quad + k_{38} + k_{40} + k_{42} + k_{43} + k_{45} + 2k_{46} + 2k_{51} + k_{52} + k_{54}, \\ n_2 &= 2k_6 + k_9 + k_{13} + k_{17} + k_{19} + k_{21} + k_{23} + k_{27} + 2k_{32} + 2k_{33} + k_{37} \\ &\quad + k_{39} + k_{41} + k_{43} + k_{44} + k_{45} + 2k_{47} + 2k_{50} + k_{53} + k_{54}. \end{aligned}$$

Now, Eqs. (6), (8) and (9) give

$$\mathbb{E} V_{T_3[3], c_3}^{2k} = 6^{2-2k} \sum_{\sum_1^{54} k_i = k} (-1)^{k'} 2^{k''} 3^{k'''} \binom{k}{k_1, \ldots, k_{54}} \prod_{i=0}^{2} \frac{l_i! \, m_i! \, n_i!}{(l_i + m_i + n_i + 3)!}.$$

Again, we list the first five even moments of the area of a random triangle in our tetrahedron $T_3$, where one point lies in $c_3$:

$$\mathbb{E}V_{T_3[3],c_3}^2 = \frac{7}{2\,400} \approx 2.91667 \cdot 10^{-3},$$

$$\mathbb{E}V_{T_3[3],c_3}^4 = \frac{11}{529\,200} \approx 2.07861 \cdot 10^{-5},$$

$$\mathbb{E}V_{T_3[3],c_3}^6 = \frac{2\,839}{10\,973\,491\,200} \approx 2.58714 \cdot 10^{-7},$$

$$\mathbb{E}V_{T_3[3],c_3}^8 = \frac{29\,419}{6\,224\,027\,040\,000} \approx 4.72668 \cdot 10^{-9},$$

$$\mathbb{E}V_{T_3[3],c_3}^{10} = \frac{4\,134\,139}{36\,352\,301\,290\,905\,600} \approx 1.13724 \cdot 10^{-10}.$$

Note that the area of a triangle in $T_3$, where one vertex is fixed to lie in $c_3$, is not larger than $\sqrt{3}/6$. Hence, we want to approximate the square root function in the interval $[0, 1/12]$ by a polynomial

$$Q(x) = \sum_{i=0}^{l} b_i x^i$$

of degree $l \in \mathbb{N}$ such that $Q(x) > \sqrt{x}$ for all $x \in [0, 1/12]$. Then we have that

$$\mathbb{E}V_{T_3[3],c_3} \leq \mathbb{E}Q(V_{T_3[3],c_3}).$$

We use Lemma 2 in [4], which yields that for $m \in \mathbb{N}$, $l = 2m + 1$ and $0 < x_0 < \ldots < x_m$, the system of equations

$$Q(x_j^2) = x_j \text{ and } Q'(x_j^2) = \frac{1}{2x_j} \quad \text{for } j = 0, \ldots, m$$

determines uniquely a polynomial $Q(x) = \sum_{i=0}^{l} b_i x^i$ with the property $Q(x) \geq \sqrt{x}$ for all $x \in \mathbb{R}$. We get a lower bound on $\mathbb{E}Q(V_{T_3[3],c_3})$ via the finite dimensional linear program

$$\min_Q \mathbb{E}Q(V_{T_3[3],c_3}) = \min_{b_i} \sum_{i=0}^{l} b_i \mathbb{E}V_{T_3[3],c_3}^{2i} \quad \text{s.t. } Q(x_\ell^2) \geq x_\ell, \ x_\ell \in \left[0, \frac{\sqrt{3}}{6}\right],$$

$$\ell = 0, \ldots, L.$$

For $l = 14$ and $L = 200$ equidistant points $x_\ell \in [0, \frac{\sqrt{3}}{6}]$, we numerically compute via Matlab and CVX that

$$\mathbb{E}Q(V_{T_3[3],c_3}) > 0.04699,$$

and this approximation is not sufficient with consideration of (7).

For $l = 15$ and $L = 1\,000$, we solve the above linear program, compute the interpolation nodes with the square root function numerically, and rationalize these points to

$$\{x_0, \ldots, x_7\} = \left\{\frac{1}{45}, \frac{1}{17}, \frac{1}{11}, \frac{1}{8}, \frac{1}{6}, \frac{1}{5}, \frac{3}{13}, \frac{7}{27}\right\}.$$

Using these points, we solve the interpolation problem from Lemma 2 in [4] and finally get a polynomial $Q_{\text{cert}}(x) = \sum_{i=0}^{15} b_i x^i$ of degree 15 with explicitly given rational coefficients $b_0, \ldots, b_{15}$ and the property $Q_{\text{cert}} \geq \sqrt{x}$. Inserting the even moments computed above, we see that

$$\mathbb{E} V_{T_3[3],c_3} \leq \mathbb{E} Q_{\text{cert}}(V_{T_3[3],c_3})$$
$$= \tfrac{1\,364\,371\,357\,215\,003\,705\,813\,514\,810\,174\,100\,904\,243\,878\,770\,680\,297\,601\,217\,509\,788\,018\,385\,106\,185\,895\,467\,950\,649\,976\,445\,591\,316\,797\,579\,615\,557\,279\,107\,985\,989}{29\,065\,551\,851\,580\,796\,169\,800\,375\,743\,177\,957\,971\,396\,211\,070\,941\,311\,273\,271\,339\,321\,255\,098\,209\,543\,350\,425\,288\,830\,624\,996\,871\,957\,546\,778\,230\,784\,000\,000\,000\,000}$$
$$< 0.046942.$$

## 5 Proof of Theorem 2

### 5.1 Some basic considerations

In the first part of the proof, we use the results on the monotonicity of the moments of volumes of random full-dimensional simplices, which have already been stated in Thorem 1. It follows with Lemma 2 that, for $d \geq 2$ and $2 \leq n \leq d$, monotonicity of the map $K \to \mathbb{E} V_{K[n]}^k$ for $K$ ranging over all $d$-dimensional convex bodies $K$ does not hold unless $n = 2$ and $k \in \mathbb{N}$ or $n = 3$ and $k \in \{1, 2\}$.

### 5.2 Proof of Proposition 1

Now we give a counterexample for $d \geq 2$ and $n = 2$. Here we prove that random chords in a specific right-angled triangle are suitable for our purposes by giving the explicit values.

We consider the ratio $r(k) := \mathbb{E} V_{T_2[2],c_2}^k / \mathbb{E} V_{T_2[2]}^k$ for the triangle $T_2$ with vertices $(0,0), (0,1)$ and $(1,0)$ and $c_2 = (1/2, 1/2)$ the midpoint of its hypotenuse. It holds:

$$r(k) = \frac{(k+3)(k+4)}{2^{k+3} + 2^{k/2+2}}.$$

The ratio $r(k)$ is clearly monotonely decreasing, and $r(1)$ is smaller than 1. Therefore, $r(k)$ is smaller than 1 for each $k \in \mathbb{N}$ and Proposition 1 holds.

### 5.3 Proof of Proposition 2

We can settle the case $n = 3, k = 2$ by comparing the second moments of the areas of our two different types of random triangles in a tetrahedron directly. We computed in Sect. 4:

$$\mathbb{E}V^2_{T_3[3],c_3} = \frac{7}{2\,400} = 0.002916... < 0.005625 = \frac{9}{1\,600} = \mathbb{E}V^2_{T_3[3]}.$$

Consequently, with Lemmata 1 and 2, we have shown that the map $K \mapsto \mathbb{E}V^2_{K[3]}$ is not monotone if $K$ ranges over all $d$-dimensional convex bodies for $d \geq 3$.

In the same section, we have shown that

$$\mathbb{E}V_{T_3[3],c_3} < 0.046942 < \mathbb{E}V_{T_3[3]}.$$

Hence, we also have a counterexample for $k = 1$, $n = 3$ and $d \geq 3$, and the proof of Proposition 2 as well as that of Theorem 2 are completed.

## 6 Proof of Theorem 3

Let $d \geq 3$, $2 \leq n \leq d$ and $k \in \mathbb{N}$. If $d > 3$, $n \neq 3$ or $k > 2$, there exists a pair of convex bodies $K, L \subseteq \mathbb{R}^{d-1}$ such that $K \subset L$, but $\mathbb{E}V^k_{K[n]} > \mathbb{E}V^k_{L[n]}$. As in the proof of Lemma 2, for $\varepsilon > 0$, we define two $d$-dimensional convex bodies by

$$K_\varepsilon := K \times [0, \varepsilon], \quad L_\varepsilon := L \times [0, \varepsilon].$$

In the sequel, we will show that $\mathbb{E}\mathcal{V}^k_{K_\varepsilon[n]}$ and $\mathbb{E}\mathcal{V}^k_{L_\varepsilon[n]}$ converge to $\mathbb{E}V^k_{K[n]}$ and $\mathbb{E}V^k_{L[n]}$, resp., for $\varepsilon$ tending to zero. Therefore, there exists an $\varepsilon_0 > 0$ such that $\mathbb{E}\mathcal{V}^k_{K_{\varepsilon_0}[n]} > \mathbb{E}\mathcal{V}^k_{L_{\varepsilon_0}[n]}$, but $K_{\varepsilon_0} \subseteq L_{\varepsilon_0}$.

Let $X_0, \ldots, X_{n-1}$ be random points in bd $K_\varepsilon$ and

$$N_0 = \#\{X_i \in K \times \{0, \varepsilon\}\}, \quad N_1 = \#\{X_i \in \text{bd } K \times [0, \varepsilon]\}.$$

Then $N_0 + N_1 = n$. It holds for $0 \leq i \leq n - 1$:

$$\mathbb{P}(X_i \in K \times \{0, \varepsilon\}) = \frac{2 \operatorname{vol}_{d-1} K}{2 \operatorname{vol}_{d-1} K + S(K)\varepsilon},$$

$$\mathbb{P}(X_i \in \text{bd } K \times [0, \varepsilon]) = \frac{S(K)\varepsilon}{2 \operatorname{vol}_{d-1} K + S(K)\varepsilon},$$

if $S(K)$ is the surface area of the boundary of $K$. Therefore, we have that

$$\begin{aligned}\mathbb{E}\mathcal{V}^k_{K_\varepsilon[n]} &= \mathbb{P}(N_1 = 0)\,\mathbb{E}\mathcal{V}^k_{K \times \{0,\varepsilon\}[n]} + \mathbb{P}(N_1 > 0)\,\mathbb{E}_{N_1>0}\mathcal{V}^k_{K_\varepsilon[n]} \\ &= \frac{2^n \operatorname{vol}_{d-1} K^n}{(2 \operatorname{vol}_{d-1} K + S(K)\varepsilon)^n}\,\mathbb{E}\mathcal{V}^k_{K \times \{0,\varepsilon\}[n]} \\ &\quad + \left(1 - \frac{2^n \operatorname{vol}_{d-1} K^n}{(2 \operatorname{vol}_{d-1} K + S(K)\varepsilon)^n}\right) \mathbb{E}_{N_1>0}\mathcal{V}^k_{K_\varepsilon[n]}.\end{aligned} \quad (10)$$

If $\varepsilon$ tends to zero, $\mathbb{P}(N_1 > 0)$ and hence the second summand in (10) tend to zero. Now, an argument similar to that in the proof of Lemma 2 yields that $\mathbb{E}\mathcal{V}^k_{K \times \{0,\varepsilon\}[n]}$ converges to $\mathbb{E}V_{K[n]}^k$, and the same is true for $K$ replaced by $L$, which completes the proof.

**Acknowledgements** Open access funding provided by Paris Lodron University of Salzburg.